# A multilevel Analysis of Saudi Arabian Student 8th Grade Mathematics Achievement TIMSS 2011

**Maha Al-Mutairi [1]\*, Khaled Bennour [2]**

1. Department of Mathematics, Shaqra University, SA; malmataery@su.edu.sa
2. Department of Statistics and Operations Research College of Science, King Saud University, SA; kbennour@ksu.edu.as
\* Correspondence: malmataery@su.edu.sa

**Abstract:** This article uses the hierarchical linear modelling (HLM) technique to explain causally the Mathematics Performance achievement of students in Saudi Arabia. Particularly, the HLM technique was applied to the TIMSS 2011 data set where five variables (Home educational resources (BSBGHER), Like learning mathematics (BSBGSLM), Self-confidence in mathematics (BSBGSCM), Engaged in mathematics learning (BSBGEML), Value learning math (BSBGSVM)) at the student level and three variables (Emphasis on academic success (BCBGEAS), School discipline and safety (BCBGDAS), and Instruction affected by mathematics resources shortages (BCBGMRS) at the school level, were used to build the hierarchical linear model so as to predict the status of mathematically 8th grader. The final model suggested that all the student level factors are found to be significant but their impact on achievement do not vary significantly across the population of schools, i.e. BSBGHER, BSBGSLM, BSBGSCM, BSBGEML, and BSBGSVM significantly predicted the status of Mathematics Performance achievement of students in Saudi Arabia. At the school level, it is found that BCBGDAS and BCBGMRS have a significant impact on performance. However, a scale point increase in the availability of school resources for mathematics decreased achievement by 4 points.

## 1. Introduction

It is well known that mathematical performance achievement is influenced by various factors, such as educational, psychological, biographical, social, among others. These factors can be categorized as students and school variables. In this article, we consider to the TIMSS 2011 data set to investigate the mathematical performance achievement in Saudi Arabia students using five student level variables and three school level variables.

Considerable studies have been done to investigate trends in mathematics achievement and the factors effecting mathematics learning and performance—e.g., [1-2,15,21,25-26,35]. For example, [21] investigated the factors of mathematics achievement including students' gender, age, ethnicity, their family socioeconomic status and school characteristics. In [26] the effects of school, students' attitudes and their beliefs in mathematics learning on students' performance were studied. Mathematics beliefs and self-concept were also considered in [15] and [36]. While [2] studied gender differences in mathematics achievement among high-school students.

Student engagement is another important factor which is defined as the level of participation, and intrinsic interest that a student shows at the school. It relies on students' behaviour at their schools such as persistence, effort, motivation, positive learning values, enthusiasm, and interest, see [8]. Various studies have displayed that student engagement is fundamentally essential in promoting achievement. This because engaging students during the learning process leads to success and more learning, both inside and outside school (e.g., [27,33].

The literature shows a variety of studies in which gender differences in mathematics were considered, see [7,12, 18,30], just to mention a few. These studies have demonstrated that gender differences in mathematics performance can be regarded as small in many countries. Furthermore, these studies provide evidence that the magnitude of gender differences has declined compared to previous decades. [7] analyzed two data sets from the TIMSS 2003 and the PISA 2003 studies to check cross-national patterns of gender differences in mathematical achievement, attitudes, and affect and assessed the links of these patterns to gender equity at the national scale. The results of this study demonstrated that the gender gap in mathematics continues in some countries. Despite the similarities between boys' and girls' achievements, it appears that boys feel more confident and less anxious in their mathematical abilities than girls. In addition, boys are more extrinsically and intrinsically to do well in mathematics than girls, which is closely consistent with related research results in the literature (e.g., [18]). Also, boys scored one third of a standard deviation higher than girls on mathematics self-concept and self-efficacy [7].

The effect of students' SES (e.g., parents' education and home educational resources) on achievement in science were excessively investigated [5,9,22,23,31,34,37]. These studies demonstrated that students from homes where their parents have a higher level of education and have more educational resources tend to perform better in science in comparison to those students their parents have lower levels of education and have less educational resources. [6] used a multilevel modelling technique to study students' achievement. The results revealed that family background characteristics accounted for 68.33% of the total variance in students' achievement [4].

Developing a positive attitude in mathematics learning has been excessively recognized. Indeed, many countries have set it as one of the main goals of mathematics education at schools. For example, in the case of Singapore, the national mathematics curriculum states "mathematics education aims to enable pupils to develop positive attitudes towards mathematics, including confidence, enjoyment and perseverance" (Ministry of Education [MOE], 2000, p. 9). The academic emphasis of school is another essential variable in explaining student achievement. Setting high academic targets for students leads to suitable learning environment which motivates students to work hard and higher academic achievements [17]. Literature that is closely related to the relationship of academic emphasis and achievement leads to consistent all levels of education's results, i.e., elementary, middle, and high school, academic emphasis and achievement were positively related, even controlling for socioeconomic factors (see, [10,16]).

The school discipline and safety characteristics explain the variance in student achievement among schools. At schools in which the disciplinary climate is strong, students usually perform better both behaviourally and academically [20]. There are many studies addressing the influence of school safety conditions on student's achievement. Based on these studies, violence has been found to hinder cognitive, social, and emotional development [28]. In violent schools, it appears that students have less time to focus on academic activities as they are paying more attention to other factors and personal safety issues, see [3,28]. So, it can conclude that unsafe school conditions have a negative impact on students' academic achievements. The relationship between school resources, (e.g., textbooks, computers, calculators, the number of pupils per teacher) and student achievement is one of the most debated issues in education which is of particular interest to policy-makers who are responsible for making decisions regarding the allocation of resources to schools. There are inconsistent results about the relationship of school resources and academic achievement. While there are studies which concluded that there is no strong and continuous link between school resources and the academic performance of students (see,[13]), some studies showed that expenditures per student had a relatively large degree of positive effect on the academic performance of students [14]. Good attendance at school (GAS) is the evaluation of the school principals concerning the seriousness of students' behaviour which is measured by: arriving late at school, skipping class and absenteeism. The findings, internationally, show that the average achievement is higher in schools where these problems are not serious (see, [22]). A cross-national study indicates that GAS accounted for a slight portion of the variance in the achievement of Singaporean eighth-graders ([19]), indicating that these problems are not serious among the Singaporean schools. In the 2007 TIMSS, Singapore, with only 4% of schools, was one the lowest countries having serious problem with students' behavior. However, the achievement average of this group of schools was lower by 211 scale points compared to schools where the students' behaviour was not serious.

This article aims to investigate the factors that affect the mathematics performance achievement of students in Saudi Arabia by using of HLM analysis to TIMSS 2011 data [1]. HLM is a comprehensive statistical technique for analyzing hierarchical structures such as students nested within schools [21,32]. Through this approach, the factors that influenced Mathematics Performance achievement are examined from both student and school perspectives. Undoubtedly, identifying which factors influence students' academic achievement is important to educational stakeholders, especially for educational decision-makers who can make use of these findings to guide both policy and practice. The specific research questions for this study are formulated as follows:

1. How do 8th grade students' mathematics performance achievement vary between students within school and across schools?
2. What factors at the student level significantly contribute to influencing students' Mathematics Performance achievement?
3. As well as controlling for student variables, what factors at the school level significantly contribute to influencing of students' mathematics achievement Performance?

The rest of the paper is structured as follows. First, we describe the Methods and the TIMSS 2011 Saudi Arabia data, including all Student Level (BSBGHER, BSBGSLM, BSBGSCM, BSBGEML and BSBGSVM) variables and School Level ((BCBGEAS, BCBGDAS and BCBGMRS) variables. Also, the descriptive statistics of all of such variables is presented. Second, two-level HLM analysis is used to estimate the model's parameters. Finally, we state our conclusion.

## 2. Materials and Methods

This study uses the TIMSS 2011 Saudi Arabia data. Various explanatory variables at the student level and the school level are used. The outcome variable for this study is students' Mathematics Performance achievement in TIMSS 2011. TIMSS 2011 employed five plausible values to estimate the Mathematics Achievement of each student. At the student level, a total of five variables are included in the analyses namely, Home educational resources (BSBGHER), Like learning mathematics (BSBGSLM), Self-confidence in mathematics (BSBGSCM), Engaged in mathematics learning (BSBGEML), Value learning math (BSBGSVM). At the school level, a total of three variables was included in the analyses extracted from the Saudi Arabia school database. These variables are: Emphasis on academic success (BCBGEAS), School discipline and safety (BCBGDAS), Instruction affected by mathematics resources shortages (BCBGMRS).

**Table 1.** Explanations for Independent Variables at student and school levels.

| Students Level | |
|---|---|
| **Variable** | **Description** |
| Home educational resources (BSBGHER) | This scale is based on 8th-grade students' responses to the following variables: number of books in the home; educational aids in the home (computer, study desk/table for own use, dictionary); and parents' education (mother and father) [1=few resources, 2=some resources, 3=many resources]. |
| Like learning maths (BSBGSLM) | Students like learning mathematics: The scale was created by TIMSS and based on students' responses to the following five statements: a) I enjoy learning mathematics; b) I wish I did not have to study mathematics; c) Mathematics is boring; d) I learn many interesting things in mathematics; e) I like mathematics [1=don't like learning math's, 2=somewhat like learning math's, 3=like learning math's]. |
| Self-confidence in maths (BSBGSCM) | Students' confidence in mathematics: The scale was created by TIMSS and based on students' responses to the following seven statements: a) I usually do well in mathematics; b) Mathematics is harder for me than for many of my classmates; c) I am just not good at mathematics; d) I learn things quickly in mathematics; e) I am good at working out difficult mathematics problems; f) My teacher tells me I am good at mathematics; g) Mathematics is harder for me than any other subject [1=not confident, 2=somewhat confident, 3=confident]. |
| Engaged mathematics learning (BSBGEML) | Engaged mathematics learning: The scale was created by TIMSS and based on students' responses to the following five statements: a) I know what my teacher expects me to do; b) I think of things not related to the lesson (reverse coded); c) My teacher is easy to understand; d) I am interested in what my teacher says; and e) My teacher gives me interesting things to do [low=1, medium=2, high=3]. |
| Value learning math (BSBGSVM) | Students' value in mathematics: The scale was created by TIMSS and based on students' responses to the following six statements: a) I think learning mathematics will help me in my daily life; b) I need mathematics to learn other school subjects ; c) I need to do well in mathematics to get into the university of my choice ; d) I need to do well in mathematics to get the job I want ; e) I would like a job that involves using mathematics; f) It is important to do well in mathematics [value =1, Somewhat Value =2, Do Not Value=3] |
| **School level** | |
| Emphasis on academic success (BCBGEAS) | School emphasis on academic success: The index was created by TIMSS and based on students' responses to the following five statements given by school principals: a) Teachers' understanding of the school's curricular goals; b) Teachers' degree of success in implementing the school's curriculum; c) Teachers' expectations for student achievement; d) Parental support for student achievement; and e) Students' desire to do well in school [1=medium, 2=high, 3=very high]. |
| INSTRUCTION AFFECTED BY MATHEMATICS | School resources: The index was created by TIMSS and based on principals' responses related to how much |

| | | |
|---|---|
| RESOURCE SHORTAGES (BCBGMRS) | capacity is available to provide instruction affected by a shortage or inadequacy of the following statements: Instructional materials (e.g., textbooks); Supplies (e.g., papers, pencils); School buildings and grounds; Heating/cooling and lighting systems; Instructional space (e.g., classrooms); Technologically competent staff; computers for instruction; Teachers with a specialization in mathematics; Computer software for mathematics instruction; Library materials relevant to mathematics instruction; Audio-visual resources for mathematics instruction; Calculators for mathematics instruction [1= affected a lot, 2=somewhat affected, 3=not affected]. |
| Discipline and safety of school (BCDGDAS) | School discipline and safety: The index was created by TIMSS and based on students' responses to the following five statements: a) This school is located in a safe neighborhood; b) I feel safe at this school; c) This school's security policies and practices are sufficient; d) The students behave in an orderly manner; and e) The students are respectful of the teachers [1=moderate problems, 2=minor problems, 3=hardly any problems]. |

Table 2. gives descriptive statistics (Mean and Std. Deviation) of Student Level and School Level variables for Saudi Arabia in TIMSS 2011.

**Table 2.** Descriptive statistics of Student Level and School Level variables for Saudi Arabia in TIMSS 2011.

| Level | Variable | Mean | Std. Deviation |
|---|---|---|---|
| Student Level | Home educational resources | 9.35 | 1.96 |
| | Student like learning mathematics | 10.05 | 2.10 |
| | Student ENGAGED with mathematics | 10.35 | 1.90 |
| | Student value learning mathematics | 10.17 | 2.00 |
| | Student CONFIDENCE with mathematics | 10.61 | 1.98 |
| School Level | School discipline and safety | 9.68 | 2.59 |
| | School emphasis on academic success principal reports | 9.89 | 2.20 |
| | Instruction affected by mathematics resources shortages | 9.33 | 1.39 |

## 3. Data Analysis

HLM is used for two-level HLM analysis (see [29]) because of the nested structure of the data and the sample design. A model building process is applied to examine the likelihood of the selected student and school variables in influencing the students' Mathematics Performance achievement modelled for Saudi Arabia students. The model building process involved the inclusion and examination of student level variables (Level 1) and followed by testing the direct and moderating effects of school level (Level 2) variables on the criterion variable.

First, the proportions of variance of student Mathematics Achievement Performance at the student and school levels are examined (i.e., fully unconditional model—Model A). Second, student variables were added to the model as Level 1, and non-significant variables are removed (model trimming). This resultant model (Model B) answers the question to what extent are the considered student variables likely to influence Mathematics Performance achievement. Third, school variables are added to model B as Level 2 explanatory variables (Model C). The final model (Model C) answers the question, which school variables contribute to explanatory effects on Mathematics Performance achievement, after controlling for student characteristics. By comparing Model B and Model C, a good understanding of how the student-level variables operate to influence Mathematics Performance achievement within schools of different characteristics can be estimated.

## 4. HLM Analysis

The unconditional model (Model A) is actually equivalent to a one-way ANOVA model with random effects:
Level 1:
$$\text{Mathematics achievement} = y_{ij} = \beta_{0j} + r_{ij}$$
Level 2:
$$\beta_{0j} = \gamma_{00} + u_{0j}$$

where Mathematics Achievement, $y_{ij}$ is the Mathematics performance achievement for student i in school j. The variance of $r_{ij}$, the variability of random error at the student level, represents the variance of Mathematics Performance achievement between students within the school, and denoted σ2, and the variance of $u_{0j}$, denoted τ00, is the variance of Mathematics Performance achievement between schools. The γ00 is the grand mean of Mathematics Performance achievement of all students. The estimates of $\gamma_{00}$, $\sigma^2$ and $\tau_{00}$ are given in Table 3. We obtain a significant non-zero grand-mean mathematics achievement score, $\gamma_{00} = 394.3$ with se=3.53. The level-1 variance estimate showed significant mathematics achievement score variation across students within a school, σ2 =5711.6. The level-2 variance shows a significant variance in the mathematics achievement means across schools, τ00=3128. The Intra-class correlation (ICC; denoted as ϱ) was calculated for the unconditional model to explore the relative school differences. Mathematically, the ICC is defined as τ00/(τ00+ σ2). Therefore, the ICC in this study was 0.3539, meaning that the variation between schools accounted for 35% of the total variance of Mathematics Performance achievement. When we use proc mixed with SAS, we obtain $\gamma_{00} = 395.9$ with a se=4.51. As expected, the coefficient estimates differ slightly but the standard error obtained by the macro is smaller as proc mixed assumes that schools are selected by a simple random sampling.

Now we add student level predictors to explain variation in mathematics achievement and obtain model B. We use the following model:
Mathematics achievement= $y_{ij} = \beta_{0j} + \beta_{1j}(\text{BSBGHER}) + \beta_{2j}(\text{BSBGSLM}) + \beta_{3j}(\text{BSBGEML}) + \beta_{4j}(\text{BSBGSCM}) + \beta_{5j}(\text{BSBGSVM}) + r_{ij}$

where
$$\beta_{0j} = \gamma_{00} + u_{0j}$$
$$\beta_{1j} = \gamma_{10} + u_{1j}$$
$$\beta_{2j} = \gamma_{20} + u_{2j}$$
$$\beta_{3j} = \gamma_{30} + u_{3j}$$
$$\beta_{4j} = \gamma_{40} + u_{4j}$$
$$\beta_{5j} = \gamma_{50} + u_{5j}$$

with $\text{var}(u_{1j}) = \tau_{11}$, $\text{var}(u_{2j}) = \tau_{22}$, $\text{var}(u_{3j}) = \tau_{33}$, $\text{var}(u_{4j}) = \tau_{44}$, $\text{var}(u_{5j}) = \tau_{55}$.

Each school mathematics achievement is identified by six parameters: the intercept $\beta_{0j}$ and the slopes $\beta_{1j}$, $\beta_{2j}$, …,$\beta_{5j}$. As predictors are centered around grand mean, the intercept represents the grand mean. Theses parameters vary across schools and all the estimated variance of the slopes are non-significant. We are not able to reject the hypothesis that $\tau_{ii} = 0$, for $i = 1,2,…,5$. Hence the relationship between level 1 predictors and mathematics achievement within schools does not vary significantly across schools. While the fixed component for the intercept is still statistically significant, and its value (395.4) has changed very little. Also, all the five student-level variables are significant. Note that student confidence with mathematics was the strong predictor of mathematics achievement in Saudi Arabia; a one scale-point increase in confidence in mathematics increased achievement by 18.5 points. In a second position but far behind we found three predictors. For home educational resources, liking learning and value learning, a one scale-point increase increased mathematics score by less than 4 points. However, students who were interested in what teacher says or said that teacher gives interesting things to do perform less than students who have less positive attitude; a one scale-point increase decreased mathematics score by more than 2 points.

Note that by comparing the $\sigma^2$ between model A and model B, we obtain the proportion reduction in variance or "variance explained" at level 1. Recall that Proportion of variance explained is

($\sigma^2$ of model A - $\sigma^2$ of model B)/$\sigma^2$ of model A

The estimated proportion of variance between students within the school explained by the model with five student predictors is

(5711.6-4408.3)/5711.6=0.23

Which means, adding the five predictors of mathematics achievement reduced the within- school variance by 23%.

Now model C can be constructed by adding school level predictors. In this case, the mathematics achievement can be represented as

Mathematics achievement= $y_{ij} = \beta_{0j} + \beta_{1j}(\text{BSBGHER}) + \beta_{2j}(\text{BSBGSLM}) + \beta_{3j}(\text{BSBGEML}) + \beta_{4j}(\text{BSBGSCM}) + \beta_{5j}(\text{BSBGSVM}) + r_{ij}$,

Where $\beta_{0j} = \gamma_{00} + \gamma_{01}\text{BCBGEAS} + \gamma_{02}\text{BCBGDAS} + \gamma_{03}\text{BCBGMRS} + u_{0j}$, with
$$\beta_{1j} = \gamma_{10} + u_{1j}\beta_{1j} = \gamma_{10} + u_{1j}$$
$$\beta_{2j} = \gamma_{20} + u_{2j}$$
$$\beta_{3j} = \gamma_{30} + u_{3j}$$
$$\beta_{4j} = \gamma_{40} + u_{4j}$$
$$\beta_{5j} = \gamma_{50} + u_{5j}$$

The proportion reduction in variance or "variance explained" at level 2 is
$$(\tau_{00} \text{ of model } B - \tau_{00} \text{ of model } C)/\tau_{00} \text{ of model } B.$$

The estimated proportion of variance between schools explained by the model with five student predictors is (2273.6-1998.1)/2273.6=0.12. Hence 12% of the between-school variance in mathematics score is explained by the three level 2 predictors. After removing the effect of the five level 1 and the three level 2 variables, the correlation between scores in the same school is slightly reduced since $\rho = \tau_{00}/(\tau_{00} + \sigma^2) = 1998.1/(1998.1 +$ 4355.6) = 0.31. Table 3 presents parameter estimates of student and school levels of the three models A, B, and C.

**Table 3.** Parameter estimate of models A, B, and C.

| Variable | Model A | | Model B | | Model C | |
|---|---|---|---|---|---|---|
| | Estimate | SE | Estimate | SE | Estimate | SE |
| Student Level | | | | | | |
| INTERCPT $\gamma_{00}$ | 394.3*** | 3.53 | 395.4*** | 3.06 | 395.5*** | 2.95 |
| Home educational resources | | | 3.8*** | 0.77 | 3.7*** | 0.76 |
| Student like learning mathematics | | | 3.8*** | 0.87 | 4.0*** | 0.88 |
| Student ENGAGED with mathematics | | | -2.5* | 1.01 | -2.7** | 1.01 |
| Student value learning mathematics | | | -3.6*** | 0.87 | -3.8*** | 0.90 |
| Student CONFIDENCE with mathematics | | | 18.5*** | 0.86 | 18.4*** | 0.88 |
| School Level | | | | | | |
| School discipline and safety | | | | | -1.4 | 1.19 |
| School emphasis on academic success | | | | | 7.9*** | 1.44 |
| Instruction affected by mathematics resources shortages | | | | | -4.0* | 1.84 |
| Level 1 variance $\sigma^2$ | 5711.6 | | 4408.3 | | 4355.6 | |
| Level 2 variance $\tau_{00}$ | 3128.0 | | 2273.6 | | 1998.1 | |
| * $p < .05$. ** $p < .01$. *** $p < .001$. | | | | | | |

Based on the findings in Table 3, the associations between the student-level variables and Mathematics Performance achievement in Model C are quite similar to those in Model B. All the five level 1 variables are significant. The mathematics achievement still strongly affects by the student confidence in mathematics. There is no statistically significant relation between school discipline and mathematics achievement. The remaining two variables, identified as the school-level variables, have a significant impact on score. However, mathematics resources shortages are found to have a negative influence on Mathematics Performance achievement.

## 5. Conclusions

In this chapter, the hierarchical linear modelling technique is used to explain causally the Mathematics Performance achievement of students in Saudi Arabia. We found that 35% of the variance in students' Mathematics Achievement is due to the variance of performance between schools and this also corresponds to the intra-class correlation. The results of PISA 2006 show that the intra-class correlation ranges from 0.06 in Finland to 0.61 in Hungary. This ranks Saudi Arabia in the medium position for the variability in performance linked to differences between schools. All the student level factors are found to be significant but their impact on achievement do not vary significantly across the population of schools. At the school level, it is found that School discipline and safety and instruction affected by mathematics resources shortage have a significant impact on performance. However, a scale point increase in the availability of school resources for mathematics decreased achievement by 4 points.

**Funding:** This research received no specific grant from any funding.

**Data Availability Statement:** TIMSS 2011 data set which is available at

https://timssandpirls.bc.edu/timss2011/international-database.html

**Conflicts of Interest:** The authors declare no conflict of interest.